\documentclass[11pt]{article}
\usepackage{amsmath, amsfonts, amssymb}
\usepackage{graphicx, amsthm,color}

\providecommand{\U}[1]{\protect\rule{.1in}{.1in}}

\newtheorem{theorem}{Theorem}[section]
\newtheorem{corollary}[theorem]{Corollary}
\newtheorem{lemma}[theorem]{Lemma}
\newtheorem{proposition}[theorem]{Proposition}
\theoremstyle{definition}
\newtheorem{definition}[theorem]{Definition}

\newtheorem{example}[theorem]{Example}

\let\epsilon\varepsilon

\begin{document}

\title{A determination theorem in terms of the metric slope}
\author{Aris Daniilidis, David Salas}

\begin{center}
{\LARGE Determination of functions by metric slopes}
\end{center}

\medskip

\begin{center}
{\large \textsc{Aris Daniilidis, David Salas}}
\end{center}

\bigskip

\noindent\textbf{Abstract.} We show that in a metric space, any continuous
function with compact sublevel sets and finite metric slope is uniquely
determined by the slope and its critical values.

\vspace{0.55cm}

\noindent\textbf{Key words} Metric slope, critical value, determination of a function.

\vspace{0.55cm}

\noindent\textbf{AMS Subject Classification} \ \textit{Primary} 49J52 ;
\textit{Secondary} 30L15, 35K90, 37C10, 58E05

\section{Introduction}

\label{sec:intro}

Let $(M,d)$ be a metric space. The \emph{metric slope} (also known as strong
slope) of a function $f:M\rightarrow\mathbb{R}$ at a point $x\in M$ is given
by
\begin{equation}
|\nabla f|(x)=
\begin{cases}
\phantom{jaimeortega}0\,,\quad & \text{ if }x\text{ is an isolated
point,}\smallskip\\
\displaystyle\limsup_{y\rightarrow x}\frac{(f(x)-f(y))^{+}}{d(x,y)}\,,\quad &
\text{ otherwise,}
\end{cases}
\label{eq:Def-MetricSlope}
\end{equation}
where $(f(x)-f(y))^{+}=\max\{f(x)-f(y),0\}$ is the standard asymmetric norm in
$\mathbb{R}$. This notion was first introduced by De Giorgi, Marino and
Tosques in \cite{DeGiorgiMarinoTosques1980} to extend steepest
descent curves to metric spaces. Since then, the metric slope has been popularized
to study descent curves in abstract settings: %\textit{e.g.} \cite{DegiovanniMarinoTosques1985, MarinoSacconTosques1989},
see the monograph \cite{AmbrosioGigliSavare2008} and references therein, or the more recent works \cite{DrusvyatskiyIoffeLewis2015, MS2020}. The metric slope has also been widely used in nonsmooth analysis in relation with concepts like error bounds (\cite{AzeCorvellec2004} \textit{e.g.}) or metric regularity and
subdifferential calculus (see \textit{e.g.} monograph \cite{Ioffe2017} and references therein).\smallskip\newline
One of the main properties that makes the metric slope a suitable generalization from vector spaces to metric spaces is that it coincides with the norm of steepest descent vectors in both, the convex and the smooth
settings. Specifically,

\begin{itemize}
\item Whenever $X$ is a normed space and $f:X\to\mathbb{R}$ is a
differentiable function at $x\in X$, one has that $|\nabla f|(x) = \|\nabla
f(x)\|$ (see, e.g. \cite[Chapter~3.1.2]{Ioffe2017}).

\item Whenever $X$ is a Banach space and $f:X\to\mathbb{R}$ is a continuous
convex function, one has that
\[
|\nabla f|(x) = \|\partial f(x)^{\circ}\|,\quad\forall x\in X,
\]
where $\partial f(x)^{\circ}$ is the element of minimal norm of the
Moreau-Rockafellar subdifferential of $f$ at $x$ (see, e.g., \cite[Theorem~5]{FabianHenrionKrugerOutrata2010}).
\end{itemize}

Recently, an ostensively unexpected result was discovered: In Banach spaces,
the metric slope contains sufficient information to determine lower semicontinuous
convex functions that are bounded from below. This determination theorem was
first established for convex functions of class $\mathcal{C}^{2}$ in
\cite[Theorem~3.8]{Daniilidis2018GradientFlows}, then for convex functions of class
$\mathcal{C}^{1}$ for which (global) minima exist\footnote{J.-B. Baillon, Personal communication in 2018} and it has eventually
been extended to Hilbert spaces in~\cite[Corollary~3.2]{PerezArosSalasVilches2020} and to Banach spaces in~\cite[Theorem~5.1]{ThibaultZagrodny2021}. The main result can be summarized as follows:

\begin{theorem}
[slope determination -- convex case]\label{thm:salas}Let $X$ be a
Banach space and $f,g:X\rightarrow\mathbb{R}$ be two convex
continuous and bounded from below functions such that
\[
|\nabla f|(x)=|\nabla g|(x),\quad\forall x\in X.
\]
Then there exists $c\in\mathbb{R}$ such that $f=g+c$.
\end{theorem}

The main strategy to obtain the above determination result when $X$ is a Hilbert space was to control the
difference of the values of the convex functions $f,g$ at an arbitrary point
$x$ by the limiting value $\inf f-\inf g$ by means of a suitable steepest
descent curve issued from~$x$. The argument uses the fact that the only
possible critical value of a convex function is its infimum together with the
existence of subgradient descent curves issued from any point of the
space.\smallskip

The extension of Theorem~\ref{thm:salas} to Banach spaces has instead been obtained through an abstract determination result formulated in metric spaces and based on the notion of \emph{global slope} (see, e.g., \cite[Definition~1.2.4]{AmbrosioGigliSavare2008}). Concretely, it has been shown in \cite[Corollary~4.1]{ThibaultZagrodny2021} that two continuous, bounded from below functions are equal up to a constant, provided they have the same (finite) global slope at every point. The notion of global slope is difficult to handle (it is based on global information, instead of local), but for convex functions in Banach spaces it does coincide with the (local notion of) metric slope and as a consequence Theorem~\ref{thm:salas} follows.\smallskip

Focusing on the metric slope, our goal in this work is to provide minimal sufficient conditions, beyond the convex case, under which for two functions
$f,g:M\rightarrow\mathbb{R}$ over a metric space $M$ the following implication
holds:
\begin{equation}
|\nabla f|(x)=|\nabla g|(x),\quad\forall x\in M\implies f=g\text{ up to an
additive constant~?}\label{eq:DeterminationImp}
\end{equation}
Although the aforementioned strategies fail in the nonconvex case even for
bounded, (globally) Lipschitz real-analytic functions defined in $\mathbb{R}$
(see Example~\ref{ex:Arctan}), in this work we show that neither convexity nor
the lineal structure of the ambient space are necessary conditions for
establishing such a determination result of local nature. The main idea is to replace the
steepest descent curves, used in \cite{PerezArosSalasVilches2020}, by discrete descent paths, constructed via transfinite
induction over the ordinals, that lead any point of the space towards a
critical point for the slope. These latter points serve as points of
comparison for the values of the functions $f$ and $g$. To this end, we
consider functions for which the construction of such discrete paths is
possible and the existence of critical points is ensured: these requirements
are fulfilled by the class of continuous functions with finite slope and
compact sublevel sets. Quite notably, the topology $\tau$ on the space $M$
under which the function is continuous and the sublevel sets are compact does
not necessarily have to be its metric topology: it can be stronger, weaker or
even not comparable.\smallskip

Let us now fix our setting and
terminology. From now on, we denote by $M$ a metric space endowed with a
distance function $d$. For a function $f:M\rightarrow\mathbb{R}$ over $M$, we
denote by $[f\leq\alpha]$ its sublevel set at the value $\alpha\in\mathbb{R}$,
that is,
\[
\lbrack f\leq\alpha]=\{x\in M\ :\ f(x)\leq\alpha\}.
\]
We recall by \eqref{eq:Def-MetricSlope} the definition of the metric slope
$|\nabla f|$ and we denote by $\mathrm{Crit}(f)$ the set of (metric) critical
points of $f$ with respect to the slope, that is,
\[
\mathrm{Crit}(f)=\{x\in M\ :\ |\nabla f|(x)=0\}.
\]
Let further $\tau$ be any topology on $M$ (which may or may not be its metric
topology). For instance, $M$ can be a Banach space and $\tau$ its weak-topology.

\begin{definition}
[$\tau$-coercive function]A function $f:M\rightarrow\mathbb{R}$ is called
\emph{$\tau$-coercive} if the sublevel sets $[f\leq\alpha]$ are $\tau$-compact
for all $\alpha<\,\underset{x\in M}{\sup}f(x).$
\end{definition}

\noindent Notice that $\tau$-coercivity ensures the existence of global
minimizers, therefore of metric critical points. If the whole space $M$ is
itself $\tau$-compact, then every $\tau$-lower semicontinuous function is
$\tau$-sublevel compact.

\section{Main results}

\label{sec:main}

In this section we show that in any metric space $(M,d)$, any $\tau
$-continuous $\tau$-coercive function (for some topology $\tau$ on $M$) can be
determined by its slope (provided it is everywhere finite) and its critical values.

\begin{lemma}
[key lemma]\label{lemma:aris}Let $f,g:M\rightarrow\mathbb{R}$ be two
real-valued functions such that
\[
|\nabla f|(x)>|\nabla g|(x)\text{, for every }x\in M\setminus\mathrm{Crit}
(f).
\]
Then, for every $x\in M\setminus\mathrm{Crit}(f)$, there exists $z\in M$ such
that%
\[
(f-g)(x)>(f-g)(z)\qquad\text{and}\qquad f(x)>f(z).
\]

\end{lemma}

\noindent\textbf{Proof. }Let $x\in M\setminus\mathrm{Crit}(f)$ and pick
$\varepsilon>0$ such that $|\nabla f|(x)>|\nabla g|(x)+2\varepsilon.$ Then
$|\nabla f|(x)>0$ (therefore, $x$ is not an isolated point in $M$) and
$|\nabla g|(x)<+\infty.$ By (\ref{eq:Def-MetricSlope}) there exists $\delta>0$
such that
\begin{equation}
\sup_{y\in B(x,\delta)}\frac{\left(  g(x)-g(y)\right)  ^{+}}{d(x,y)}<|\nabla
g|(x)+\varepsilon. \label{aa1}%
\end{equation}
If $|\nabla f|(x)<+\infty$ there exists $z\in B(x,\delta)$ such that
\[
\frac{\left(  f(x)-f(z)\right)  ^{+}}{d(x,z)}\,>\,|\nabla f|(x)-\varepsilon.
\]
Since $|\nabla f|(x)-\varepsilon>\,|\nabla g|(x)+\varepsilon\,>0,$ it follows
that
\[
f(x)-f(z)=\left(  f(x)-f(z)\right)  ^{+}\,>\,d(x,z)\,\left(  |\nabla
g|(x)+\varepsilon\right)  .
\]
The above inequality holds true also if $|\nabla f|(x)=+\infty$ and yields
$f(x)-f(z)>0.$ Combining with \eqref{aa1} we obtain
\[
f(x)-f(z)>d(x,z)\,\left(  |\nabla g|(x)+\varepsilon\right)  >\left(
g(x)-g(z)\right)  ^{+}\geq g(x)-g(z)
\]
and the statement of the lemma follows.\hfill$\Box$

\bigskip

\begin{proposition}[strict comparison]
\label{prop:comparison} Let $f,g:M\rightarrow\mathbb{R}$ be $\tau$-continuous,
$\tau$-coercive functions (for some topology $\tau$ in $M$) and
\[
|\nabla f|(x)\,>\,|\nabla g|(x),\,\text{ for every }x\in M\setminus\mathrm{Crit}(f).
\]
Then, $\mathrm{Crit}(f)\neq\emptyset$ and for every $x\in M\setminus
\mathrm{Crit}(f)$
\begin{equation}
f(x)-g(x)\,\,>\,\,m(x):=\underset{z\in\lbrack f\leq f(x)]\cap\mathrm{Crit}(f)}{\inf
}(f-g)(z). \label{eq:aris}
\end{equation}

\end{proposition}

\noindent\textbf{Proof. }Apply Lemma~\ref{lemma:aris} to obtain $x_{0}\in M$
such that $(f-g)(x)>(f-g)(x_{0})$ and $f(x)>f(x_{0}).$ If $x_{0}
\in\mathrm{Crit}(f),$ then $(f-g)(x)>(f-g)(x_{0})\geq m(x)$ and
\eqref{eq:aris} holds. If $x_{0}\in M\diagdown\mathrm{Crit}(f),$ then we can
apply again Lemma~\ref{lemma:aris} and obtain $x_{1}\in M$ such that
$(f-g)(x_{0})>(f-g)(x_{1})$ and $f(x_{0})>f(x_{1}).$ As long as we do not meet
a critical point of $f$ we build a generalized sequence $\{x_{\alpha
}\}_{\alpha}\subset\lbrack f\leq f(x_{0})]$ over the ordinals as follows: for
every ordinal $\lambda$ for which $\{x_{a}\}_{a<\lambda}\subset M\setminus
\mathrm{Crit}(f)$ is defined:

\begin{enumerate}
\item[$(i)$.] If $\lambda=\beta+1$ is a sucessor ordinal and $x_{\beta}\in
M\setminus\mathrm{Crit}(f)$ has been defined, then we apply
Lemma~\ref{lemma:aris} for $x=x_{\beta}$ and set $x_{\beta+1}:=z$. Then
\begin{equation}
(f-g)(x_{\beta})>(f-g)(x_{\beta+1})\qquad\text{and\qquad}f(x_{\beta
})>f(x_{\beta+1}). \label{eq:suc}
\end{equation}

\item[$(ii)$.] Assume now that $\lambda$ is a limit ordinal and $\{x_{\beta}
\}_{\beta<\lambda}\subset\lbrack f\leq f(x_{0})]$ has been defined such that
\eqref{eq:suc} holds true. Since $[f\leq f(x_{0})]$ is $\tau$-compact, the
set
\[
A:=\bigcap_{\beta<\lambda}\overline{\{x_{\alpha}\ :\ \beta\leq\alpha
<\lambda\}}^{\tau}
\]
is nonempty. Pick any $x_{\lambda}\in A$. Then, for every $\beta<\lambda$, by
$\tau$-continuity of the function $f-g$ we have
\begin{align*}
(f-g)(x_{\beta})  &  =\sup\left\{  (f-g)(x)\ :\ x\in\overline{\{x_{\alpha
}\ :\ \beta\leq\alpha<\lambda\}}^{\tau}\right\} \\
&  \geq(f-g)(x_{\lambda}).
\end{align*}

\end{enumerate}

\noindent The above construction will necessarily end up to a critical
point of $f$, that is, we eventually obtain $x_{\lambda}\in\mathrm{Crit}(f)$ for some
ordinal $\lambda$. Indeed, if this does not happen, then reaching ordinals of
arbitrary cardinality (in particular, bigger than $|M|$) we deduce the
existence of two ordinals $\beta<\alpha$ for which $x_{\beta}=x_{\alpha}.$
Then our construction yields
\[
(f-g)(x_{\beta})>(f-g)(x_{\beta+1})\geq(f-g)(x_{\alpha}),
\]
which is obviously a contradiction. Therefore, $\mathrm{Crit}(f)\neq
\emptyset.$ Moreover, $\{x_{a}\}_{a}\subset\lbrack f<f(x_{0})]$ and for every
$x_{a}$ we have
\[
(f-g)(x_{0})>(f-g)(x_{\alpha})\geq\underset{z\in\lbrack f\leq f(x)]\cap
\mathrm{Crit}(f)}{\inf}(f-g)(z).
\]
The proof is complete.\hfill$\Box$

\bigskip

\begin{proposition}
[Comparison Principle]\label{thm:ComparisonPrinciple} Let $f,g:M\rightarrow
\mathbb{R}$ be $\tau$-continuous, $\tau$-coercive functions. If

\begin{enumerate}
\item[(i)] $|\nabla g|(x)\leq|\nabla f|(x)<+\infty$ for all $x\in M$, and

\item[(ii)] There exists $c\in\mathbb{R}$ such that $g(x)-f(x)\leq c$ for all
$x\in\mathrm{Crit}(f)$,
\end{enumerate}

then $g \leq f+c$.
\end{proposition}

\noindent\textbf{Proof.} The conclusion is trivial if $\mathrm{Crit}(f)=M.$ If
$\mathrm{Crit}(f)\neq M,$ set $f_{\varepsilon}=(1+\varepsilon)f$ and notice
that $f_{\varepsilon}$ is $\tau$-continuous and $\tau$-coercive. Then, for
every $x\in M$ we know that $|\nabla f_{\varepsilon}|(x)=(1+\varepsilon
)|\nabla f|(x)$ (see, e.g. \cite[Proposition 3.3]{Ioffe2017}). Thus, we have
that $\mathrm{Crit}(f_{\varepsilon})=\mathrm{Crit}(f)$ and $|\nabla
g|(x)=|\nabla f|(x)<|\nabla f_{\varepsilon}|(x)$ for every $x\in
M\setminus\mathrm{Crit}(f_{\varepsilon})$. By
Proposition~\ref{prop:comparison}, we have
\begin{align*}
g(x)  &  <f_{\varepsilon}(x)-\underset{z\in\lbrack f_{\varepsilon}\leq
f_{\varepsilon}(x)]\cap\mathrm{Crit}(f)}{\inf}(f_{\varepsilon}-g)(z)\\
&  =f_{\varepsilon}(x)+c-\underset{z\in\lbrack f_{\varepsilon}\leq
f_{\varepsilon}(x)]\cap\mathrm{Crit}(f)}{\inf}\varepsilon f(z)\\
&  =f(x)+\varepsilon\left[  \underbrace{ f(x)-\underset{z\in\lbrack
f_{\varepsilon}\leq f_{\varepsilon}(x)]\cap\mathrm{Crit}(f)}{\inf}
f(z)}_{\text{uniformly bounded for all } \varepsilon}\right]  +c.
\end{align*}
By taking $\varepsilon\searrow0$, we conclude that $g(x)\leq f(x)+c$. The
proof is complete. \hfill$\Box$

\bigskip

By applying twice Proposition~\ref{thm:ComparisonPrinciple}, we can deduce the
main result of this work.

\begin{theorem}
[Determination Theorem]\label{thm:DeterminationTheorem} Let $M$ be a metric
space and $\tau$ be a topology over $M$. Let $f,g:M\rightarrow\mathbb{R}$ be
two $\tau$-continuous and $\tau$-coercive functions such that

\begin{enumerate}
\item[(i)] $|\nabla f|(x)=|\nabla g|(x)<+\infty$ for every $x\in M$.

\item[(ii)] There exists $c\in\mathbb{R}$ such that $g(x)-f(x)=c$ for all
$x\in\mathrm{Crit}(f)$.
\end{enumerate}

Then, $g=f+c$.
\end{theorem}

A direct albeit important corollary is that, whenever $M$ is a compact metric
space, the determination theorem applies for continuous functions with finite
slope, by taking $\tau$ to be the topology defined by the metric $d$.

\begin{corollary}
\label{cor:DeterminationTheoremCompact} Let $M$ be a compact metric space and
let $f,g:M\to\mathbb{R}$ be two continuous functions with finite slope such that

\begin{enumerate}
\item[(i)] $|\nabla f|(x) = |\nabla g|(x)$ for all $x\in M$.

\item[(ii)] There exists $c\in\mathbb{R}$ such that $g(x) - f(x) = c$ for all
$x\in\mathrm{Crit}(f)$.
\end{enumerate}

Then, $g = f+c$.
\end{corollary}

The distinction between the metric~$d$ for which the slope is being computed,
and the topology~$\tau$ over which continuity and sublevel coercivity are considered, yields 
the following corollary over (possibly infinite-dimensional) reflexive Banach spaces.

\begin{corollary}
\label{cor:DeterminationBanach} Let $X$ be a reflexive Banach space and let
$f,g:X\to\mathbb{R}$ be two weak-continuous and weak-sublevel coercive
functions such that

\begin{enumerate}
\item[(i)] $|\nabla f|(x) < +\infty$ for every $x\in M$.

\item[(ii)] $|\nabla f|(x) = |\nabla g|(x)$ for all $x\in M$.

\item[(iii)] There exists $c\in\mathbb{R}$ such that $g(x) - f(x) = c$ for all
$x\in\mathrm{Crit}(f)$.
\end{enumerate}

Then, $g=f+c$.
\end{corollary}

\section{Pertinence of the assumptions and illustrative examples}

\label{sec:Counterexamples} The slope, by its own, does not
provide enough information to verify~\eqref{eq:DeterminationImp} for arbitrary
continuous functions, unless additional assumptions are imposed. In this
section we illustrate the relevance of our assumptions by means of
counterexamples.\smallskip
\newline Let us first observe that even in the framework of
Theorem~\ref{thm:salas} (convex case), the assumption that the convex
functions are bounded from below is pertinent: indeed, without this
assumption, one could obtain an effortless counterexample by considering two
different linear functionals (elements of the dual space) of the same norm.
The simplest concrete example is to consider the 1--dimensional case
$X=\mathbb{R}$ and the (linear) functions $f(t)=t$ and $g(t)=-t$ for
all $t\in\mathbb{R}$.\smallskip\newline 
Concurrently, the following example
shows that mere boundedness from below without convexity is not enough, even
for real-analytic, bounded functions. This example illustrates, in case of
absence of convexity, the need to assume existence of global minima (or at
least, of critical points). We recall that in the framework of Theorem~\ref{thm:DeterminationTheorem}, 
this is ensured by the coercivity assumption.

\begin{example}
\label{ex:Arctan} Let us consider the functions $f,g:\mathbb{R}\rightarrow
\mathbb{R}$ given by
\[
f(t)=\arctan(t)\quad\mbox{ and }\quad g(t)=-\arctan(t)
\]
\begin{figure}[h]
\centering
\includegraphics[scale=0.2]{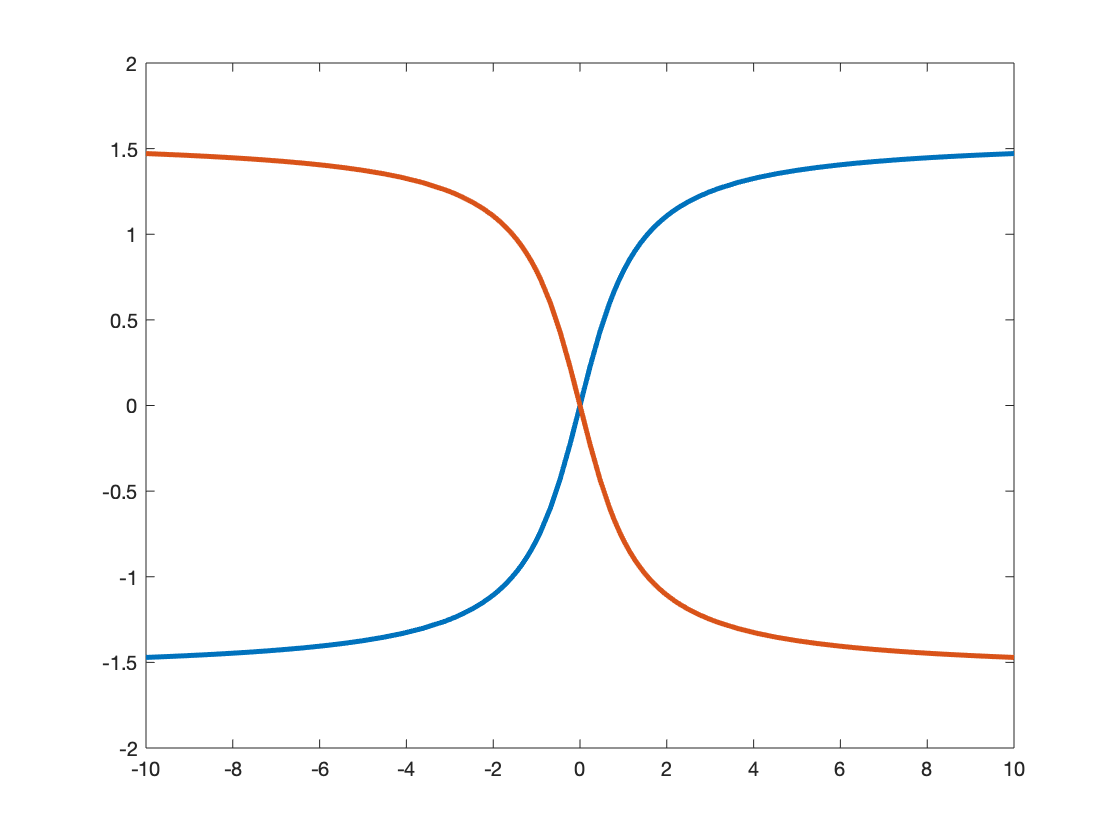}
\caption{The functions $f(t) = \arctan(t)$ (in blue) and $g(t) = -\arctan(t)$
(in red).}
\label{fig:arctan}
\end{figure}
Since both functions are continuously differentiable, we have that
\[
|\nabla f|(t)=|f^{\prime}(t)|=\frac{1}{t^{2}+1}=|g^{\prime}(t)|=|\nabla
g|(t).
\]
Nevertheless, the desired implication \eqref{eq:DeterminationImp} does not
hold (see Figure \ref{fig:arctan}).\qed
\end{example}

The following example shows that for (nonconvex) smooth coercive functions in
$\mathbb{R}$, implication \eqref{eq:DeterminationImp} might fail if we lack
information on the behavior of the functions on the set of critical values.

\begin{example}
Consider the function $f:\mathbb{R}\rightarrow\mathbb{R}$ given by
\begin{equation}
f(x)=
\begin{cases}
(x+\pi/2)^{2}-1,\qquad & \text{ if }x<-\pi/2\\
\sin(x),\qquad & \text{ if }x\in\lbrack-\pi/2,\pi/2]\\
(x-\pi/2)^{2}+1,\qquad & \text{ if }x>\pi/2
\end{cases}
\label{eq:FunctionSquareAndSinus}
\end{equation}
It is not hard to see that $f$ is of class $\mathcal{C}^{1}$ and coercive,
with
\[
f^{\prime}(x)=
\begin{cases}
2(x+\frac{\pi}{2}), & \text{ if }x<-\frac{\pi}{2}\\
\cos(x), & \text{ if }x\in\lbrack-\frac{\pi}{2},\frac{\pi}{2}]\\
2(x-\frac{\pi}{2}), & \text{ if }x>\frac{\pi}{2}
\end{cases}
\quad\text{and\quad}|\nabla f|(x)=
\begin{cases}
\cos(x), & \text{ if }x\in\lbrack-\frac{\pi}{2},\frac{\pi}{2}]\smallskip\\
2(|x|-\frac{\pi}{2}), & \text{ if }x\notin\lbrack-\frac{\pi}{2},\frac{\pi}{2}]
\end{cases}
\]

\begin{figure}[h]
\centering
\includegraphics[scale=0.2]{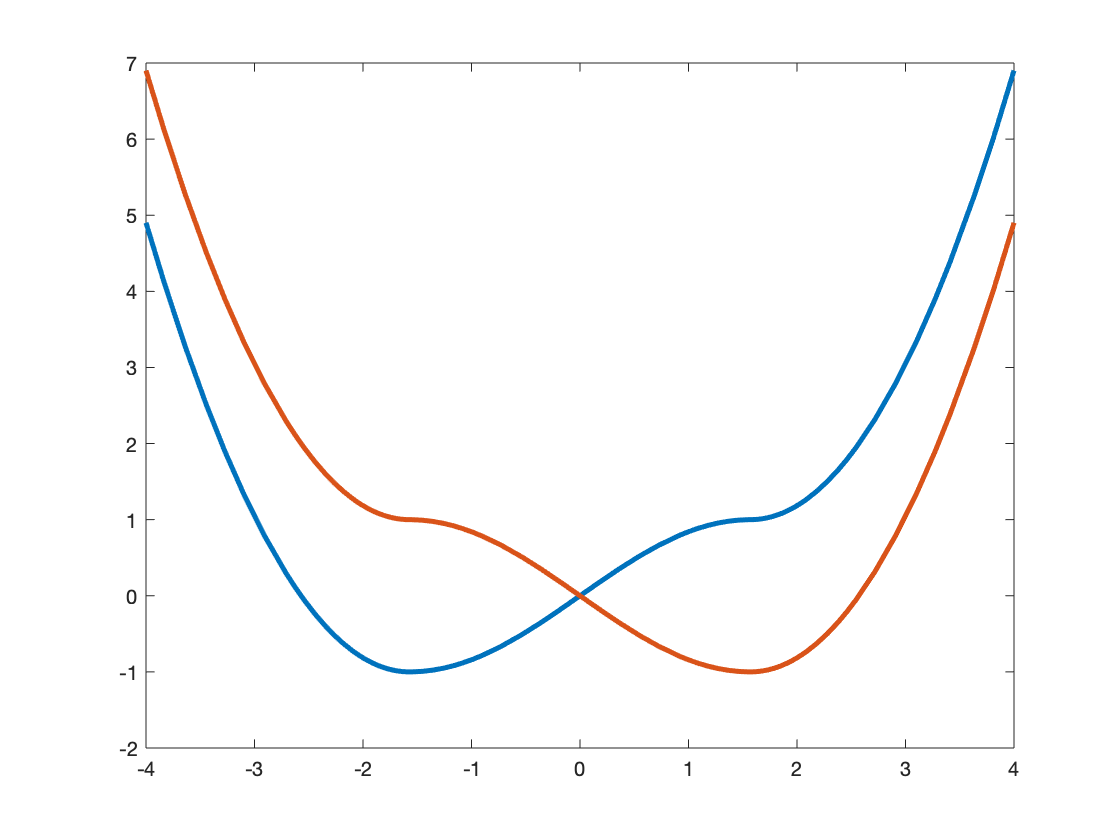}
\caption{The functions $f(x)$ given by \eqref{eq:FunctionSquareAndSinus} (in
blue) and $g(x) = f(-x)$ (in red).}
\label{fig:SquareAndSinus}
\end{figure}
\noindent By considering $g:\mathbb{R}\rightarrow\mathbb{R}$ given by $g(x)=f(-x)$, we
get that $|\nabla f|(x)=|\nabla g|(x)$ for every $x\in\mathbb{R}$, but
\eqref{eq:DeterminationImp} fails, as illustrates
Figure~\ref{fig:SquareAndSinus}.\qed

\end{example}

The last example illustrates the problem of allowing the slope to take the
value $+\infty$ quite often. In this example, the functions are continuous,
have compact sublevel sets and the same critical values, yet
\eqref{eq:DeterminationImp} fails.

\begin{example}
\label{ex:CantorStaircase} Let us consider the classical Cantor Staircase
function $\mathfrak{c}:[0,1]\rightarrow\lbrack0,1]$, which has directional
derivatives equal to either $0$ or $+\infty.$

\begin{figure}[h]
\centering
\includegraphics[scale=0.2]{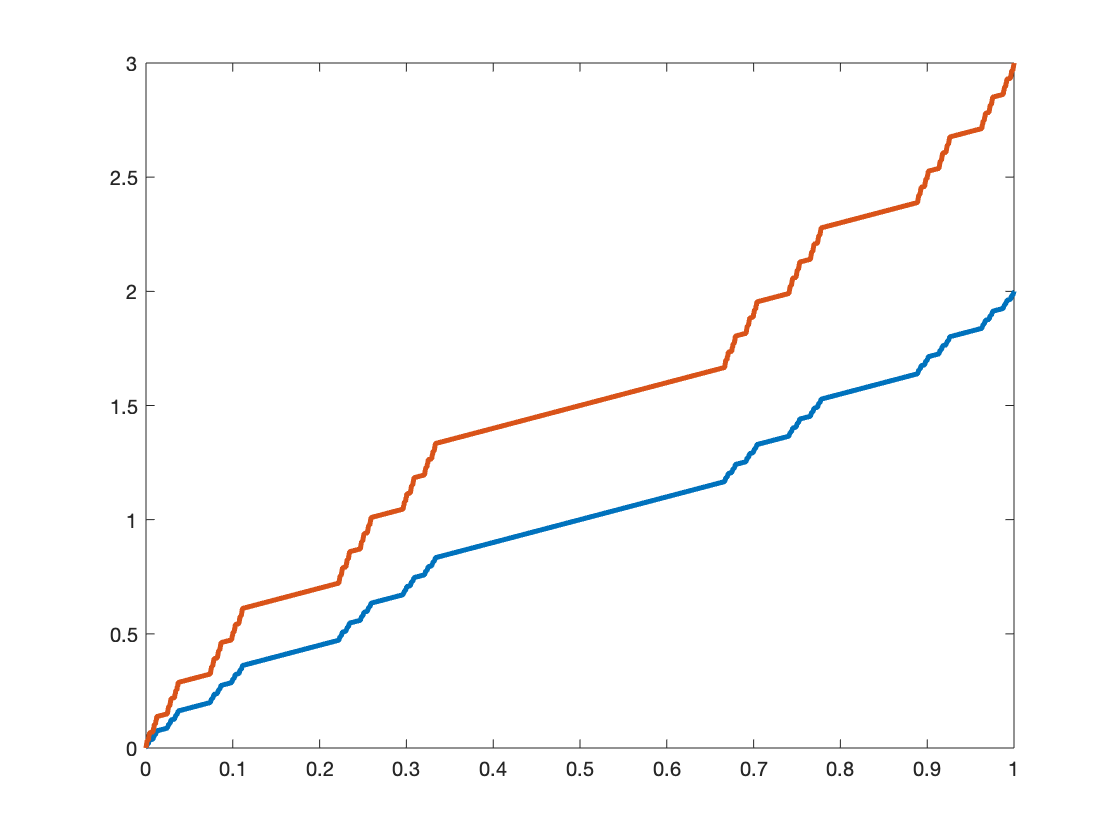}
\caption{The functions $f(t) = \mathfrak{c}(t) + t$ (in blue) and $g(t) = 2
\mathfrak{c}(t) +t$ (in red).}
\label{fig:CantorStaircase}
\end{figure}
Define $f,g:[0,1]\rightarrow\mathbb{R}$ by
\[
f(t)=\mathfrak{c}(t)+t\qquad\text{ and }\qquad g(t)=2\mathfrak{c}(t)+t.
\]
Both functions are continuous and coercive and it is not hard to see that
\[
|\nabla f|(t)=|\nabla g|(t)\in\{1,+\infty\},\qquad\text{for all }t\in(0,1]
\]
and
\[
|\nabla f|(0)=|\nabla g|(0)=0.
\]
Therefore, $\mathrm{Crit}(f)=\mathrm{Crit}(g)=\{0\}$ and $f(0)=g(0)=0$.
Nevertheless, the determination theorem \eqref{eq:DeterminationImp} fails to
hold, precisely because the slope is not finite everywhere.\qed

\end{example}

\bigskip

\noindent\textit{Final conclusion}. In this work, we established a determination result
for functions in a general framework, based on three complementary hypotheses:
Finite slope (at each point), continuity and sublevel compactness for a
suitable topology and knowledge of the critical values. Our main result shows
that in this case the metric slope contains sufficient first-order information
to determine the function. The first hypothesis seems to be necessary, since
if the slope can take infinite values, one cannot control the variation of the
functions over those points, as Example~\ref{ex:CantorStaircase} illustrates.
The third hypothesis acts as a boundary condition: since the slope does not
provide information over the set of critical points, in order to verify
\eqref{eq:DeterminationImp}, it is necessary to impose that the
functions are equal (up to an additive constant) over this set, which is
nonempty thanks to the second assumption. Nonewithstanding, it seems there is
room for improvement (weakening) of this second assumption. A perspective of
this work is to provide an alternative condition over spaces where compactness
is not present, like general Banach spaces.

\bibliographystyle{plain}
\bibliography{Biblio.bib}

\noindent\rule[0pt]{5cm}{1pt}

\vspace{0.5cm}

\noindent Aris DANIILIDIS

\medskip

\noindent DIM--CMM, UMI CNRS 2807\newline Beauchef 851, FCFM, Universidad de
Chile \smallskip

\noindent E-mail: \texttt{arisd@dim.uchile.cl} \newline\noindent
\texttt{http://www.dim.uchile.cl/\symbol{126}arisd}

\medskip

\noindent Research supported by the grants: \newline CMM AFB170001, FONDECYT
1211217 (Chile), ECOS-ANID C18E04 (Chile, France).\newline\vspace{0.5cm}

\noindent David SALAS

\medskip

\noindent Instituto de Ciencias de la Ingenieria, Universidad de
O'Higgins\newline Av. Libertador Bernardo O'Higgins 611, Rancagua, Chile
\smallskip

\noindent E-mail: \texttt{david.salas@uoh.cl} \newline\noindent
\texttt{https://www.researchgate.net/profile/David-Salas-8} \medskip

\noindent Research supported by the grant: \newline CMM AFB170001, FONDECYT
3190229 (Chile)

\end{document}